\documentclass[twoside]{article}

\usepackage{imsart,amsfonts}

\newtheorem{theorem}{Theorem}

\newtheorem{definition}{Definition}

\newtheorem{lemma}[theorem]{Lemma}
\newtheorem{proposition}[theorem]{Proposition}

\def\bds{\begin{displaystyle}}
\def\eds{\end{displaystyle}}

\setattribute{infoline}   {text} {}

\begin{document}

\begin{frontmatter}

\title{A Family of non-Gaussian Martingales with Gaussian Marginals}
\runtitle{Martingales with Gaussian Marginals}
\author{Kais Hamza}
\address{School of Mathematical Sciences, Monash University}
\author{Fima C. Klebaner}
\address{School of Mathematical Sciences, Monash University}
\runauthor{K. Hamza and F.C. Klebaner}

\begin{abstract}
We construct a family of non-Gaussian martingales the marginals of which are all Gaussian.
We give the predictable quadratic variation of these processes and show they do not have continuous
paths. These processes are Markovian and inhomogeneous in time, and we give their infinitesimal generators.
Within this family we find a class of piecewise deterministic pure jump processes and describe the laws of
jumps and times between the jumps.
\end{abstract}

\begin{keyword}[class=AMS]
\kwd[Primary ]{60G44}
\kwd{60J25}
\kwd{60J75}
\kwd[; secondary ]{91B70}
\end{keyword}

\begin{keyword}
\kwd{Martingales with Gaussian marginals}
\kwd{Time-inhomogeneous Markov processes}
\end{keyword}

\end{frontmatter}

\section{Introduction}

This paper is concerned with constructing an alternative process
to the Brownian motion; a martingale with Gaussian marginals, yet
not a Gaussian process. The construction of martingales with given marginals
has significance in financial modelling (see for example \cite{BibSkoSor05}, \cite{Cam04},
\cite{CarMad05} and \cite{MadYor02}).
In \cite{HamKle06}, the authors investigate the existence of an alternative model for
which the Black-Scholes formula holds true. The case of the
Bachelier formula raises the question of the existence of
a non-Gaussian martingale with Gaussian marginals.

Our solution to this problem is based on the elementary observation that for
$Y$ and $\xi$ independent standard Gaussian random variables, the
distribution of $Z=\sqrt rY+\sqrt{1-r}\xi$ is a standard Gaussian
random variable for any value  of $r\in [0,1]$. This allows us to
randomize $r$ and construct a family of Markovian martingales with
Gaussian marginals.

The question of constructing (Markovian) martingales with given
marginals has seen considerable interest in recent years, mostly
initiated by the paper of Madan and Yor (\cite{MadYor02}).

In \cite{MadYor02}, the authors give three different approaches: a continuous
martingale, a time-changed Brownian motion, and a cosntruction
that uses Az\`ema-Yor's solution to the Skorokhod embedding
problem. These constructions are applied to a number of special
cases including the case with Gaussian marginals. But out of the
three constructions, only the Skorokhod embedding approach yields
a non-Gaussian martingale. The other two reduce to a construction
of a Brownian motion. The continuous martingale approach looks for
a process of the type
$$X_t = \int_0^t\sigma(X_s,s)dW_s,$$
for which the marginal densities $g(x,t)$ are $N(0,t)$. Writing
the forward equation for these densities,   it follows that
$\sigma^2\equiv1$ (see \cite{MadYor02}), and  $X_t$ is a Brownian
motion. In the time change approach
$$X_t = B_{L_t}$$
where $L$ is an increasing process -- in fact $L$ is assumed
to be an increasing Markov process with inhomogeneous independent
increments, independent of the Brownian motion $B$. The
assumption of Gaussian marginals implies
$$\mathbb{E}\big[e^{i\lambda X_t}\big] =
\mathbb{E}\big[e^{-(\lambda^2/2)L_t}\big]=e^{-(\lambda^2/2)t},$$
and it follows that $L_t=t$. But the Skorokhod embedding approach
of Madan and Yor yields a discontinuous and time-inhomogeneous
Markov process.

Our approach is different to the above. It uses basic principles,
and has the advantage of producing an entire family of processes
indexed by an (infinite) family of subordinators. The construction
produces a family discontinuous and time-inhomogeneous Markov processes.
We obtain the quadratic variation of these processes and infinitesimal
generators in some cases. The richness of the family has the
potential to allow for the imposition of specifications other than
the marginal distributions.

Note that our method can be extended to include other types of
marginal distributions, but for clarity of
presentation we choose to focus solely on the Gaussian case.

Finally, all existing approaches yield discontinuous processes
(barring the Brownian motion itself), and the question of the
existence of a non-Gaussian
continuous martingale with Gaussian marginals remains open.

\section{A Family of non-Gaussian Martingales with Gaussian
Marginals} \label{construction}

In this section we construct a (non-Gaussian) Markov martingale
$X_t$ the marginals of which are Gaussian with mean zero and
variance $t$. The existence of such process is guaranteed by a
Theorem of Kellerer (see \cite{Kel72} and \cite{MadYor02}) that
only requires the targeted marginal densities, $g(x,t)$, be
increasing in the convex order
($\mathbb{E}[f(X_t)]\geq\mathbb{E}[f(X_s)]$ for $s<t$ and $f$
convex), and have means that do not depend on $t$.

As eluded to in the introduction, the main idea of the proposed
construction is the fact that, for any triple $(R,Y,\xi)$ of
independent random variables such that $R$ takes values in
$[0,1]$, $\xi$ is standard Gaussian and $Y$ is Gaussian with mean
zero and variance $\alpha^2$, the random variable $Z =
\sigma(\sqrt{R}Y+\alpha\sqrt{1-R}\xi)$ is Gaussian with mean zero
and variance $\sigma^2\alpha^2$. However, the unconditional joint
distribution of $(Y,Z)$ is not bivariate Gaussian, as can be
verified by calculating the fourth conditional moment of $Z$ given
$Y=0$. In fact, $(Y,Z)$ is a bivariate Gaussian pair if and only
if $R$ is non-random. The martingale property of the two-step
process $(Y,Z)$ holds if and only if
$$Y=\mathbb{E}[Z|Y] = \mathbb{E}[\sigma(\sqrt{R}Y+\alpha\sqrt{1-R}\xi)|Y]
= \sigma\mathbb{E}[\sqrt{R}]Y,$$ in other words,
\begin{equation}
\mathbb{E}\big[\sqrt{R}\big] = \frac1\sigma. \label{sqrtcond}
\end{equation}   Furthermore, the
conditional distribution of $Z$ given $Y$ is
$$F_{Z|Y=y}(dz) =
\mathbb{P}[R=1]\varepsilon_{\sigma y}(dz) +
\mathbb{E}\left[\phi\left(\sigma\sqrt{R}y,\alpha^2\sigma^2(1-R),z\right)1_{R<1}\right]dz,$$
where $\varepsilon_x$ is the Dirac measure at $x$ and
$\phi(\mu,\sigma^2,\cdot)$ denotes the density of the Gaussian
distribution with mean $\mu$ and variance $\sigma^2$.

This construction of a two-step process can be extended to that of
a continuous time Markov process. Indeed, let $R_{s,t}$ be a
family of random variables indexed by $0<s\leq t$. We assume that
$R_{s,t}$ takes values in $[0,1]$, has distribution that depends
on $(s,t)$ only through $\sqrt{t/s}$, and has moment of order
$1/2$ equal to $\sqrt{s/t}$. We denote the distribution of
$R_{s,t}$ by $G_{\sqrt{t/s}}(dr)$. We shall also need a family
$\xi_{s,t}$ of standard Gaussian random variables.

The process $X_t$ is constructed as a Markov process with the
following almost sure representation of $X_t$ in terms of $X_s$,
$s<t$,
\begin{equation}
X_t = \sqrt{\frac ts}\left(\sqrt{R_{s,t}}X_s + \sqrt
s\sqrt{1-R_{s,t}}\xi_{s,t}\right).\label{asconstruct}
\end{equation}
We assume given the usual set-up of a probability space endowed
with a filtration ${\cal F}_t$ to which $X_t$ is adapted. In the
representation (\ref{asconstruct}), $R_{s,t}$ and $\xi_{s,t}$ are
assumed to be independent of each other, ${\cal F}_t$-measurable
and independent of ${\cal F}_s$.

In the sequel, we shall often write $\alpha$ for $\sqrt s$,
$\sigma$ for $\sqrt{t/s}$, $\tau$ for $\sqrt{u/t}$ and whenever
independence between the variables involved need not be
emphasized, $R_\sigma$ for $R_{s,t}$ and $R_\tau$ for $R_{t,u}$.

\begin{definition}
The family $(G_\sigma)_{\sigma\geq1}$ is a log-convolution
semi-group if the the distribution of the product of any two
independent random variables with distributions $G_\sigma$ and
$G_\tau$, is $G_{\sigma\tau}$.
\end{definition}

Define, for $\sigma\geq1$ and $R_\sigma$ distributed as $G_\sigma$,
$U_\sigma = -\ln R_\sigma$, and, for $p\geq0$, $V_p
= U_{e^p}$. If $K_p$ denotes the distribution of $V_p$, then
$(G_\sigma)_{\sigma\geq1}$ is a log-convolution semi-group if and
only if $(K_p)_{p\geq0}$ is a convolution semi-group:
$$K_0=\varepsilon_0,\quad K_p\ast K_q=K_{p+q}.$$

\begin{proposition}
Define, $P_{s,t}(x,dy)$ as,
\begin{equation}
P_{0,t}(x,dy) = \frac1{\sqrt{2\pi}\sqrt
t}\exp\left(-\frac{(y-x)^2}{2t}\right)dy,\label{D0}
\end{equation}
and for $s>0$,
\begin{eqnarray}
\lefteqn{P_{s,t}(x,dy)}\label{D1}\\
& = & \gamma(\sigma)\varepsilon_{\sigma x}(dy) +
\int_{[0,1)}\frac1{\sqrt{2\pi}\sqrt t\sqrt{1-r}}
\exp\left(-\frac{(y-\sigma\sqrt
rx)^2}{2t(1-r)}\right)G_{\sigma}(dr)dy\nonumber\\
& = & \gamma(\sigma)\varepsilon_{\sigma x}(dy) +
\mathbb{E}\left[\phi\left(\sigma\sqrt{R_\sigma}x,\alpha^2\sigma^2(1-R_\sigma),y\right)1_{R_\sigma<1}\right]dy,\nonumber
\end{eqnarray}
where $R_\sigma$ is distributed as $G_\sigma$ and $\gamma(\sigma)=G_\sigma(\{1\})$.

If $(G_\sigma)_{\sigma\geq1}$ is a log-convolution semi-group then
for any $u>t>s>0$ and any $x$,
\begin{equation}
\int P_{s,t}(x,dy)P_{t,u}(y,dz) = P_{s,u}(x,dz)\label{CK1}
\end{equation}
and, for any $u>t>0$,
\begin{equation}
\int P_{0,t}(0,dy)P_{t,u}(y,dz) = P_{0,u}(0,dz).\label{CK0}
\end{equation}
\label{CKpos}
\end{proposition}
{\bf Proof:} We prove (\ref{CK1}) and (\ref{CK0}) by showing that
the almost sure formulation (\ref{asconstruct}) is consistent.
\begin{eqnarray*}
X_u & = & \tau\left(\sqrt{R_{t,u}}X_t + \sigma\alpha
\sqrt{1-R_{t,u}}\xi_{t,u}\right)\\
& = & \tau\left(\sqrt{R_{t,u}}\sigma\left(\sqrt{R_{s,t}}X_s +
\alpha \sqrt{1-R_{s,t}}\xi_{s,t}\right) + \sigma\alpha
\sqrt{1-R_{t,u}}\xi_{t,u}\right)\\
& = & \sigma\tau\left(\sqrt{R_{s,t}R_{t,u}}X_s + \alpha\left(
\sqrt{(1-R_{s,t})R_{t,u}}\xi_{s,t} + \sqrt{1-R_{t,u}}\xi_{t,u}\right)\right)\\
\end{eqnarray*}

Now, letting $R_{s,u}=R_{s,t}R_{t,u}$ and
$$\xi_{s,u} = \left(\frac{\sqrt{(1-R_{s,t})R_{t,u}}}{\sqrt{1-R_{s,u}}}1_{R_{s,u}<1} + 1_{R_{s,u}=1}\right)\xi_{s,t} +
\frac{\sqrt{1-R_{t,u}}}{\sqrt{1-R_{s,u}}}1_{R_{s,u}<1}\xi_{t,u},$$
we see that
$$X_u = \tau\sigma\left(\sqrt{R_{s,u}}X_s + \alpha
\sqrt{1-R_{s,u}}\xi_{s,u}\right)$$ with $R_{s,u}$ distributed as
$G_\tau(dr)$. Also, the unconditional distribution of $\xi_{s,u}$
as well as its conditional distribution given $R_{s,t}$ and
$R_{t,u}$ are standard Gaussian. This in turn implies that
$\xi_{s,u}$ is independent of $R_{s,u}$.\hfill$\Box$

\begin{proposition}[L\'evy-Khinchin Theorem]
Assume that the family $(G_\sigma)_{\sigma\geq1}$ is a
log-convolution semi-group and let $(R_\sigma)_{\sigma\geq1}$ be
independent random variables with laws $(G_\sigma)_{\sigma\geq1}$.
Let $L_\sigma(\lambda) = \mathbb{E}\left[e^{\lambda\ln
R_\sigma}\right] = \mathbb{E}\left[(R_\sigma)^\lambda\right]$ be
the moment generating function of the (positive) random variable
$U_\sigma = -\ln R_\sigma$.

For any $\sigma\geq1$, $U_\sigma = -\ln R_\sigma$ is infinitely
divisible, and
\begin{equation}
\ln L_\sigma(\lambda) = -\left[\beta\lambda + \int_0^\infty\left(1
- e^{-\lambda x}\right)\nu(dx)\right]\ln\sigma
\end{equation}
where the L\'evy measure $\nu(dx)$ satisfies $\nu(\{0\})=0$ and
$\bds\int_0^\infty(1\wedge x)\nu(dx)<\infty\eds$.
\end{proposition}

In what follows, we denote by $\psi$ the Laplace exponent of the
log-convolution semi-group $(G_\sigma)_{\sigma\geq1}$:
$$\psi(\lambda) = \beta\lambda + \int_0^\infty\left(1 - e^{-\lambda x}\right)\nu(dx).$$
The following theorem   follows form the above and the
Chapman-Kolmogorov existence result.
\begin{theorem}
Assume that the family $(G_\sigma)_{\sigma\geq1}$ is a
log-convolution semi-group with Laplace exponent
$$\psi(\lambda) =
\beta\lambda + \int_0^\infty\left(1 - e^{-\lambda
x}\right)\nu(dx).$$ If $\psi(1/2)=1$, then there exists a Markov
martingale $X_t$ starting at zero with transition probabilities
$P_{s,t}(x,dy)$ given by (\ref{CK1}) and (\ref{CK0}) the marginal
distributions of which are Gaussian with mean zero and variance
$t$.\label{existence}
\end{theorem}

\section{Path properties}

\begin{theorem}
The process $X_t$ is continuous in probability:
$$\forall c>0,\ \lim_{s\rightarrow t}\mathbb{P}[|X_t-X_s|>c] = 0.$$
\end{theorem}
{\bf Proof:} Using Lemma \ref{secondmoment} below, we write,
$$\mathbb{P}[|X_t-X_s|>c] \leq \frac1{c^2}\mathbf{E}[(X_t-X_s)^2] =
\frac1{c^2}[t - t^{1-\delta}s^\delta + t^{1-\delta}s^{\delta} - s]
= \frac{t-s}{c^2}.$$ \hfill$\Box$
\begin{lemma}
Let $\delta=\psi(1)/2$ so that $L_\sigma(1)=\sigma^{-2\delta}$.
Then
$$\mathbb{E}[(X_t-X_s)^2|X_s] = t - t^{1-\delta}s^\delta +
t^{1-\delta}s^{-1+\delta}X_s^2 - X_s^2.$$\label{secondmoment}
\end{lemma}
{\bf Proof:} Using representation (\ref{asconstruct}), we see that
\begin{eqnarray*}
\mathbf{E}\left[(X_t-X_s)^2|X_s\right]
& = & \mathbf{E}[\mathbf{E}[(X_t-X_s)^2|X_s,R_\sigma]|X_s]\\
& = & \alpha^2\sigma^2\mathbf{E}\left[1-R_\sigma\right]+
\mathbf{E}\left[(\sigma\sqrt{R_\sigma}-1)^2\right]X_s^2\\
& = & \alpha^2\sigma^2\left(1-L_\sigma(1)\right)+
\left(\sigma^2\mathbf{E}[R_\sigma]-1\right)X_s^2\\
& = & \alpha^2\sigma^2\left(1-L_\sigma(1)\right)+
\left(\sigma^2L_\sigma(1)-1\right)X_s^2\\
& = & \alpha^2\sigma^2 - \alpha^2\sigma^{2-2\delta} +
\sigma^{2-2\delta}X_s^2 - X_s^2\\
\end{eqnarray*}
\hfill$\Box$

\begin{theorem}
The (predictable) quadratic variation of $X_t$ is
$$\left<X,X\right>_t = \delta t+(1-\delta)\int_0^t\frac{X_s^2}{s}ds,$$
where $\delta = \psi(1)/2$. Furthermore, it can be obtained as a
limit in probability,
$$\left<X,X\right>_t =
\mathbb{P}\lim\sum_{k=0}^{n-1}\mathbf{E}\left[(X_{t_{k+1}}-X_{t_k})^2|X_{t_k}\right]$$
where $t_0<t_1<\ldots<t_n$ is a subdivision of $[0,t]$.
\label{quadratic}
\end{theorem}
{\bf Proof:} First note that $X_t$ is a square integrable
martingale on any finite interval $[0,T]$. In fact
$\bds\sup_{t\leq T}\mathbf{E}[X_t^2] = T\eds$. Also,
$$\mathbf{E}[X_t^2|{\mathcal{F}}_s]
= \mathbf{E}\left[t(1-R_\sigma)+\sigma^2R_\sigma X_s^2|X_s\right]
= t(1-L_\sigma(1)) + \sigma^2L_\sigma(1)X_s^2.$$ Since
$L_\sigma(1)=\sigma^{-2\delta} = s^\delta t^{-\delta}$, we find
$$\mathbf{E}[X_t^2|{\mathcal{F}}_s]
= t - t^{1-\delta}s^\delta + t^{1-\delta}s^{-1+\delta}X_s^2.$$ It
follows that
\begin{eqnarray*}
\lefteqn{\mathbf{E}\left[\left.(1-\delta)\int_0^t\frac{X_u^2}{u}du\right|{\mathcal{F}}_s\right]}\\
& = & (1-\delta)\int_0^s\frac{X_u^2}{u}du
+ (1-\delta)\int_s^t(1 - u^{-\delta}s^\delta + u^{-\delta}s^{-1+\delta}X_s^2)du\\
& = & (1-\delta)\int_0^s\frac{X_u^2}{u}du
+ (1-\delta)(t-s) - s^\delta(1 - s^{-1}X_s^2)(t^{1-\delta}-s^{1-\delta})\\
& = & (1-\delta)\int_0^s\frac{X_u^2}{u}du
+ (1-\delta)(t-s) - s^\delta t^{1-\delta} + s + t^{1-\delta}s^{-1+\delta}X_s^2 - X_s^2\\
\end{eqnarray*}
and
\begin{eqnarray*}
\lefteqn{\mathbf{E}\left[\left.X_t^2 - \delta t - (1-\delta)\int_0^t\frac{X_u^2}{u}du\right|{\mathcal{F}}_s\right]}\\
& = & t - t^{1-\delta}s^\delta + t^{1-\delta}s^{-1+\delta}X_s^2 -
\delta t
- (1-\delta)\int_0^s\frac{X_u^2}{u}du - (1-\delta)(t-s) +\\
& & s^\delta t^{1-\delta} - s - t^{1-\delta}s^{-1+\delta}X_s^2 + X_s^2\\
& = & t - \delta t - (1-\delta)\int_0^s\frac{X_u^2}{u}du - (1-\delta)(t-s) - s + X_s^2\\
& = & X_s^2 - \delta s - (1-\delta)\int_0^s\frac{X_u^2}{u}du.
\end{eqnarray*}
\hfill$\Box$

The next result states that the only continuous process that can
be constructed in the the way described in Section
\ref{construction} is the Brownian motion.

\begin{theorem}
If $R_\sigma$ is not degenerate ($R_\sigma\not\equiv\sigma^{-2}$), $X_t$
is not continuous.
\end{theorem}
{\bf Proof:} We proceed by contradiction and assume that $X_t$ is
continuous. It\^o's formula for $e^{i\lambda X_t}$ gives
\begin{equation}
e^{i\lambda X_t} = 1 + M_t - \frac{\lambda^2}2\int_0^te^{i\lambda
X_s}d\left<X,X\right>_s,\label{Ito}
\end{equation}
where $\bds M_t = \int_0^ti\lambda e^{i\lambda X_s}dX_s\eds$ is a
true martingale. In fact
\begin{eqnarray*}
\mathbf{E}\left[\left|\left<M,M\right>_t\right|\right] & = &
\mathbf{E}\left[\left|-\int_0^t\lambda^2e^{i2\lambda X_s}d\left<X,X\right>_s\right|\right]\\
& = & \mathbf{E}\left[\left|-\int_0^t\lambda^2e^{i2\lambda X_s}
\left(\delta ds+(1-\delta)\frac{X_s^2}{s}ds\right)\right|\right]\\
& \leq & \delta\lambda^2t + (1-\delta)\lambda^2\int_0^t\frac{\mathbf{E}[X_s^2]}{s}ds\\
& = & \delta\lambda^2t + (1-\delta)\lambda^2t\\
& = & \lambda^2t,
\end{eqnarray*}
since $X_s$ is $\mathbf{N}(0,s)$ and $\mathbf{E}[X_s^2]=s$.

Taking expectations in (\ref{Ito}), we obtain that
$\bds\theta(\lambda,t) = \mathbf{E}[e^{i\lambda X_t}] =
e^{-\lambda^2t/2}\eds$ must satisfy
\begin{eqnarray*}
\theta(\lambda,t) & = & 1 -
\frac{\lambda^2}2\left[\delta\int_0^t\theta(\lambda,s)ds +
(1-\delta)\int_0^t\mathbf{E}\left[X_s^2e^{i\lambda X_s}\right]ds\right]\\
& = & 1 - \frac{\lambda^2}2\left[\delta\int_0^t\theta(\lambda,s)ds
-
(1-\delta)\int_0^t\frac{\partial^2\theta}{\partial\lambda^2}(\lambda,s)ds\right].
\end{eqnarray*}
Differentiating in $t$, we get that $\theta(\lambda,t)$ must
satisfy
$$-\frac{\lambda^2}2\theta(\lambda,t) = -\frac{\lambda^2}2\left[\delta\theta(\lambda,t)
-1(1-\delta)\frac{\partial^2\theta}{\partial\lambda^2}(\lambda,t)\right],$$
that is,
$$-\frac{\lambda^2}2 = -\frac{\lambda^2}2\left[\delta
-(1-\delta)t(\lambda^2t-1)\right].$$ This, of course, can only
occur if $\delta=1$, which corresponds to
$L_\sigma(1)=\sigma^{-2}$ and $R_\sigma$ being non-random equal to
$\sigma^{-2}$. \hfill$\Box$

\section{Explicit Constructions}

Before we engage in the explicit construction of the processes outlined
in the previous sections, let us observe that these fall into one of two subclasses
according to whether or not $G_\sigma(\{1\})$ is nil, uniformly in $\sigma$.

Indeed,
$$\gamma(\sigma) = G_\sigma(\{1\}) = \lim_{\lambda\uparrow\infty}L_\sigma(\lambda)
= \lim_{\lambda\uparrow\infty}\exp\left(-\psi(\lambda)\ln\sigma\right)$$
and
$$\gamma(\sigma) = 0 \Leftrightarrow \lim_{\lambda\uparrow\infty}\psi(\lambda) = +\infty.$$

\subsection{The Case $\gamma(\sigma)>0$}

In this section we apply our construction to the case where
$\gamma(\sigma)=G_\sigma(\{1\})>0$. The processes thus obtained
are piecewise deterministic pure jump process in the sense that
between any two consecutive jumps, the process behaves according
to a deterministic function. Examples of such processes include
the case where $G_\sigma$ is an inverse log-Poisson distribution.

The interpretation of these processes as piecewise deterministic
pure jump processes requires the computation of the infinitesimal
generator.

\begin{proposition}
Let $G_\sigma$ be a log-convolution semi-group for which
$\gamma(\sigma)=G_\sigma(\{1\})>0$, $\gamma$ is differentiable at
$1$ and $\bds\lim_{\lambda\downarrow0}\psi(\lambda)=0\eds$. Then
the infinitesimal generator of $X_t$ on the set of
$C^2_0$-functions is given by
\begin{eqnarray*}
\lefteqn{A_0f(x) = \frac12f''(x)\mbox{ and for }s>0,}\\
\lefteqn{A_sf(x) = \frac{x}{2s}f'(x)}\\
& & + \frac{-\gamma'(1)}{2s}\int[
f(x+z)-f(x)]\int_{[0,1)}\phi((\sqrt r-1)x,s(1-r),z)\bar G(dr)dz,
\end{eqnarray*}
where
$$\bar G(dr) = \lim_{\sigma\downarrow1}\frac{G_\sigma(dr\cap[0,1))}{G_\sigma([0,1))}$$
is a probability measure on $[0,1)$, and the limit is understood
in the weak sense. \label{discrete}
\end{proposition}
Thus the process $X$ starts off as a Browwnian motion and,
when in $x$ at time $s$, drifts at the rate
of $x/(2s)$, and jumps at the rate of $-\gamma'(1)/(2s)$. The size
of the jump from $x$ has density $\int_{[0,1)}\phi((\sqrt
r-1)x,s(1-r),z)\bar G(dr)$, the mean of which is
$\int_{[0,1)}(\sqrt r-1)\bar G(dr)x$. In other words, while in
positive territory, $X_t$ continuously drifts upwards and has
jumps that tend to be negative. In negative region, the reverse
occurs; $X_t$ drifts downwards and has (on average) positive
jumps.

{\bf Proof:} First note that the conditional moment generating
function of $U_\sigma$ given $U_\sigma>0$ is
$$L_\sigma^*(\lambda) = \frac{L_\sigma(\lambda)-\gamma(\sigma)}{1-\gamma(\sigma)}$$
and converges to
$$\lim_{\sigma\downarrow1}L_\sigma^*(\lambda) = 1+\frac{\psi(\lambda)}{\gamma'(1)}.$$
By the (Laplace) continuity theorem, if
$\lim_{\lambda\downarrow0}\psi(\lambda)=0$ then there is exists a
probability measure on $[0,1)$, $\bar G(dr)$, such that
$$\bar G(dr) = \lim_{\sigma\downarrow1}\frac{G_\sigma(dr\cap[0,1))}{G_\sigma([0,1))}.$$

Next,
\begin{eqnarray*}
\lefteqn{\frac1{t-s}\left(\mathbb{E}[f(X_t)|X_s=x]-f(x)\right)}\\
& = & \frac1s\left[\frac{f(\sigma
x)\gamma(\sigma)-f(x)}{\sigma^2-1}\right.\\
& & \left.+ \frac1{\sigma^2-1}\int
f(y)\int_{[0,1)}\phi(\sigma\sqrt
rx,t(1-r),y)G_\sigma(dr)dy\right]\\
& = & \frac1s\left[\frac{f(\sigma
x)\gamma(\sigma)-f(x)}{\sigma^2-1}\right.\\
& & \left.+ \frac{1-\gamma(\sigma)}{\sigma^2-1}\int
f(y)\int_{[0,1)}\phi(\sigma\sqrt
rx,t(1-r),y)\frac{G_\sigma(dr)}{1-\gamma(\sigma)}dy\right]\\
\end{eqnarray*}
Letting $\sigma$ decrease to 1, we see that
\begin{eqnarray*}
\lefteqn{A_sf(x)}\\
& = & \lim_{t\downarrow
s}\frac1{t-s}\left(\mathbb{E}[f(X_t)|X_s=x]-f(x)\right)\\
& = & \frac1s\left[\frac{xf'(x)+\gamma'(1)f(x)}2 -
\frac{\gamma'(1)}2\int f(y)\int_{[0,1)}\phi(\sqrt rx,s(1-r),y)\bar
G(dr)dy\right]\\
& = & \frac{x}{2s}f'(x) + \frac{-\gamma'(1)}{2s}\int[
f(y)-f(x)]\int_{[0,1)}\phi(\sqrt rx,s(1-r),y)\bar G(dr)dy\\
& = & \frac{x}{2s}f'(x)\\
& & + \frac{-\gamma'(1)}{2s}\int[
f(x+z)-f(x)]\int_{[0,1)}\phi((\sqrt r-1)x,s(1-r),z)\bar G(dr)dz.
\end{eqnarray*}
\hfill$\Box$

Note that the domain of $A_s$ can be extended to include functions
that do not vanish at infinity, such as $f(x)=x^2$. Indeed by
Theorem \ref{quadratic}, $g_s(x) = \delta + (1-\delta)\frac{x^2}s$
solves the martingale problem for $f(x)=x^2$.

The next proposition immediately follows from the observation that
the process $X$ does not jump between times $s$ and $t$ if and only
if $X_u = \sqrt{\frac us}X_s$ for $u\in(s,t)$.
\begin{proposition}
Let $T_s$ denote the first jump time after $s>0$. Then, for any
$t>s$,
$$\mathbb{P}[T_s>t] = \gamma(\sigma),$$
where as before, $\sigma = \sqrt{t/s}$.
\end{proposition}

\subsection{The Poisson Case: $\gamma(\sigma) = \sigma^{-c}$}
In this case $\beta=0$, $\nu(dx) = c\delta_1(dx)$ with $\bds
c=\frac1{1-e^{-1/2}}\eds$, and $\bds\psi(\lambda) =
c(1-e^{-\lambda})\eds$. In other words $U_\sigma = -\ln R_\sigma$
has a Poisson distribution with mean $c\ln\sigma$.

The assumptions of Proposition \ref{discrete} are clearly
satisfied with $\bds\gamma(\sigma) = \sigma^{-c}\eds$,
$\gamma'(1)=-c$, $\bds\lim_{\sigma\downarrow1}L_\sigma^*(\lambda)
= e^{-\lambda}\eds$ and $\bar G(dr)=\varepsilon_{e^{-1}}(dr)$, so
that $X_t$ has infinitesimal generator
$$A_sf(x) = \frac{x}{2s}f'(x) + \frac{c}{2s}\int[
f(x+z)-f(x)]\phi(-x/c,s(1-e^{-1}),z)dz.$$ It jumps at the rate of
$\bds\frac{c}{2s}\eds$ with a size distributed as a Gaussian
random variable with mean $\bds-\frac{x}{c}\eds$ and variance
$s(1-e^{-1})$. The graph below shows a simulation of a path of
such a process.

\includegraphics{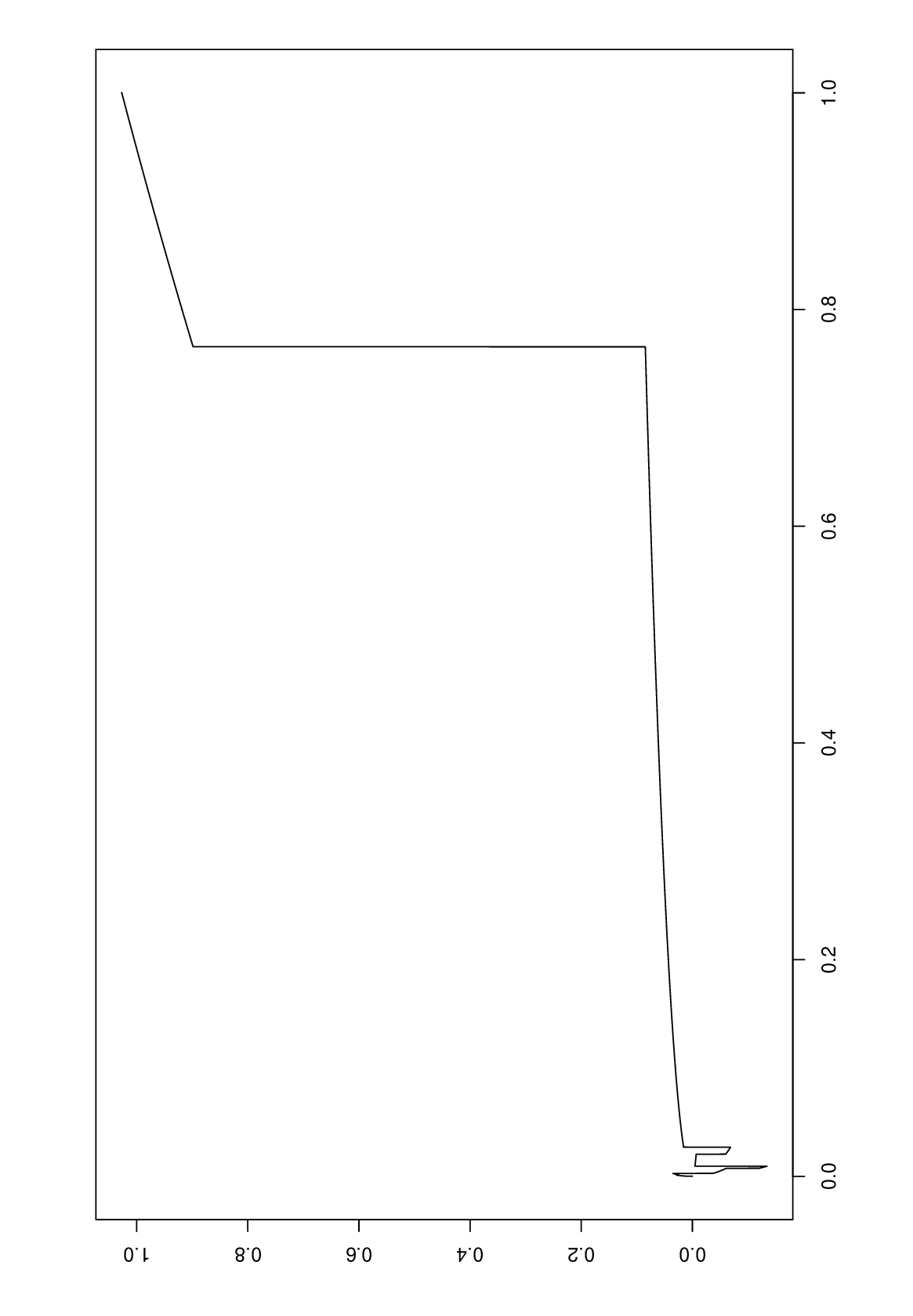}

\vspace*{6cm}

Furthermore, the law of the first jump time after $s$ is given by
$$\mathbb{P}[T_s>t] = \gamma(\sigma) = \frac{s^{c/2}}{t^{c/2}}.$$
In other words, $T_s$ is Pareto distributed (with location
parameter $s$ and scale parameter $c/2\sim1.27$). In particular,
$$\mathbb{E}[T_s] = \frac{cs}{c-2}\mbox{ and } \mathbb{E}[T_s^2] = \infty.$$

\subsection{The Case $\gamma(\sigma) = 0$}

When $\gamma(\sigma)=0$, we are only able to compute the
infinitesimal generator for functions of a  specific type.
Examples of such functions include polynomials.

\begin{proposition}
\label{generalcase}
Assume that $\beta=0$ so that
$$\psi(\lambda) = \int_0^\infty\left(1 - e^{-\lambda x}\right)\nu(dx).$$
Let $f$ be a $C^1$-function with the following property: there exist a function $N_f$
and a (signed) finite measure $M_f$ such that
$$f\left(\sigma e^{-u/2}x+\sqrt t\sqrt{1-e^{-u}}z\right) =
N_f(\sigma)\int_0^\infty e^{-\lambda u}M_f(s,x,z,d\lambda),\quad u>0,$$
and
$$\lim_{\sigma\downarrow1}N_f(\sigma) = 1.$$
Then, for any $s>0$,
\begin{eqnarray*}
\lefteqn{A_sf(x) = \frac{x}{2s}f'(x)}\\
& & + \frac1{2s}\int[
f(x+y)-f(x)]\int_0^{+\infty}\phi((e^{-\omega/2}-1)x,s(1-e^{-\omega}),y)\nu(d\omega)dy.
\end{eqnarray*}
\end{proposition}
{\bf Proof:} Let
$$C_\sigma f(u) = C_\sigma f(s,x,z,u) = f\left(\sigma e^{-u/2}x+\sqrt t\sqrt{1-e^{-u}}z\right).$$
Then, since $\gamma(\sigma)=0$, $U_\sigma$ is almost surely strictly positive and,
\begin{eqnarray*}
\lefteqn{\frac1{t-s}\left(\mathbb{E}[f(X_t)|X_s=x]-f(x)\right)}\\
& = &
\frac1s\frac1{\sigma^2-1}\int\left(\mathbb{E}\left[f\left(\sigma
e^{-U_\sigma/2}x+\sqrt
t\sqrt{1-e^{-U_\sigma}}z\right)\right]-f(x)\right)\phi(z)dz.\\
& = & \frac1s\frac1{\sigma^2-1}\left\{
\int\big(\mathbb{E}\left[C_\sigma f(U_\sigma)\right]-C_\sigma
f(0)\big)\phi(z)dz
+ \big(f(\sigma x)-f(x)\big)\right\}\\
& = & \frac1s\frac1{\sigma^2-1}\left\{
\int\mathbb{E}\left[N_f(\sigma)\int_0^\infty\left(e^{-\lambda
U_\sigma}-1\right)M_f(d\lambda)\right]\phi(z)dz
+ \left(f(\sigma x)-f(x)\right)\right\}\\
& = & \frac1s\frac1{\sigma^2-1}\left\{
N_f(\sigma)\int\int_0^\infty\left(e^{-\psi(\lambda)\ln\sigma}-1\right)M_f(d\lambda)\phi(z)dz
+ \left(f(\sigma x)-f(x)\right)\right\}\\
& = &
\frac1s\frac1{\sigma+1}\left\{\frac{N_f(\sigma)\ln\sigma}{\sigma-1}
\int\int_0^\infty\frac{e^{-\psi(\lambda)\ln\sigma}-1}{\ln\sigma}M_f(d\lambda)\phi(z)dz
+ \frac{f(\sigma x)-f(x)}{\sigma-1}\right\}.\\
\end{eqnarray*}
Taking the limit as $\sigma\downarrow1$ (that is $t\downarrow s$),
we get
$$A_sf(x) = \frac{x}{2s}f'(x) - \frac1{2s}\int\int_0^\infty\psi(\lambda)M_f(s,x,z,d\lambda)\phi(z)dz.$$
Since
$$\psi(\lambda) = \int_0^\infty\left(1 - e^{-\lambda\omega}\right)\nu(d\omega),$$
\begin{eqnarray*}
\lefteqn{A_sf(x)}\\
& = & \frac{x}{2s}f'(x) -
\frac1{2s}\int\int_0^\infty\left[\int_0^\infty\left(1 - e^{-\lambda\omega}\right)\nu(d\omega)\right]M_f(s,x,z,d\lambda)\phi(z)dz\\
& = & \frac{x}{2s}f'(x)
- \frac1{2s}\int\int_0^\infty\left[\int_0^\infty\left(1 - e^{-\lambda\omega}\right)M_f(s,x,z,d\lambda)\right]\nu(d\omega)\phi(z)dz\\
& = & \frac{x}{2s}f'(x)
+ \frac1{2s}\int\int_0^\infty\left[f\left(e^{-\omega/2}x+\sqrt s\sqrt{1-e^{-\omega}}z\right)-f(x)\right]\nu(d\omega)\phi(z)dz\\
\end{eqnarray*}
and the proof is completed by a change of variables in $z$.
\hfill$\Box$

\begin{lemma}
Let $f(x)=x^n$, then
$$f\left(\sigma e^{-u/2}x+\sqrt t\sqrt{1-e^{-u}}z\right) =
\sigma^n\int_0^\infty e^{-\lambda u}M_f(s,x,z,d\lambda)$$ where
$$M_f(s,x,z,d\lambda) = \sum_{k=0}^n\sum_{j=0}^{n-k}\frac{n!}{k!j!(n-k-j)!}(-1)^jx^ks^{(n-k)/2}z^{n-k}
\left(\varepsilon_{k/2}*m_j\right)(d\lambda)$$ and $m_j(d\lambda)$
is the $j$-order convolution of the probability measure
$$m(d\lambda) = \frac{1}{2\sqrt\pi}\sum_{n=1}^{+\infty}\frac{\Gamma(n-1/2)}{n!}\varepsilon_n(d\lambda).$$
\end{lemma}
{\bf Proof:} First, write the Taylor series of the (analytic on
$(0,1)$) function $1-\sqrt{1-x}$,
$$1-\sqrt{1-x} = \frac12\sum_{n=1}^{+\infty}\frac{\Gamma(n-1/2)}{n!\Gamma(1/2)}x^n.$$
It immediately follows that,
$$1-\sqrt{1-e^{-u}} = \frac12\sum_{n=1}^{+\infty}\frac{\Gamma(n-1/2)}{n!\Gamma(1/2)}e^{-nu}
= \int_0^\infty e^{-\lambda u}m(d\lambda),$$ where $\bds
m(d\lambda) =
\frac1{2\sqrt\pi}\sum_{n=1}^{+\infty}\frac{\Gamma(n-1/2)}{n!}\varepsilon_{n}(d\lambda)\eds$
is a probability measure. Now,
\begin{eqnarray*}
\lefteqn{f\left(\sigma e^{-u/2}x+\sqrt t\sqrt{1-e^{-u}}z\right)}\\
& = & \sigma^n\sum_{k=0}^n{n\choose k}e^{-ku/2}x^ks^{(n-k)/2}(1-e^{-u})^{(n-k)/2}z^{n-k}\\
& = & \sigma^n\sum_{k=0}^n{n\choose k}e^{-ku/2}x^ks^{(n-k)/2}\left[1-\left(1-\sqrt{1-e^{-u}}\right)\right]^{n-k}z^{n-k}\\
& = &
\sigma^n\sum_{k=0}^n\sum_{j=0}^{n-k}\frac{n!}{k!j!(n-k-j)!}(-1)^jx^ks^{(n-k)/2}z^{n-k}
e^{-ku/2}\left(1-\sqrt{1-e^{-u}}\right)^j\\
\end{eqnarray*}
The proof is ended by observing that
$$e^{-ku/2}\left(1-\sqrt{1-e^{-u}}\right)^j = \int_0^\infty e^{-\lambda u}
\left(\varepsilon_{k/2}*m_j\right)(d\lambda).$$ \hfill$\Box$

The following theorem is now proven.
\begin{theorem}
Assume that $\beta=0$. For any polynomial $f$ and any $s>0$,
\begin{eqnarray}
\lefteqn{A_sf(x) = \frac{x}{2s}f'(x)}\nonumber\\
& & + \frac1{2s}\int[
f(x+y)-f(x)]\int_0^{+\infty}\phi((e^{-\omega/2}-1)x,s(1-e^{-\omega}),y)\nu(d\omega)dy.
\label{generator}
\end{eqnarray}
\end{theorem}

\subsection{The Gamma Case: $\gamma(\sigma) = 0$}

Here $\beta=0$, $\nu(dx) = ax^{-1}e^{-bx}dx$ with $\bds
a=\frac1{\ln\left(1+\frac1{2b}\right)}\eds$ and $\bds\psi(\lambda)
= a\ln\left(1+\frac\lambda b\right)\eds$; that is
$U_\sigma$ has a gamma distribution with density
$$h_\sigma(u)  = \frac{b^{a\ln\sigma}}{\Gamma(a\ln\sigma)}u^{a\ln\sigma-1}e^{-bu}, \quad u>0,$$
and $R_\sigma$ has an inverse log-gamma distribution with density
$$g_\sigma(r)  = \frac{b^{a\ln\sigma}}{\Gamma(a\ln\sigma)}(-\ln
r)^{a\ln\sigma-1}r^{b-1}, \quad 0<r<1.$$

\includegraphics{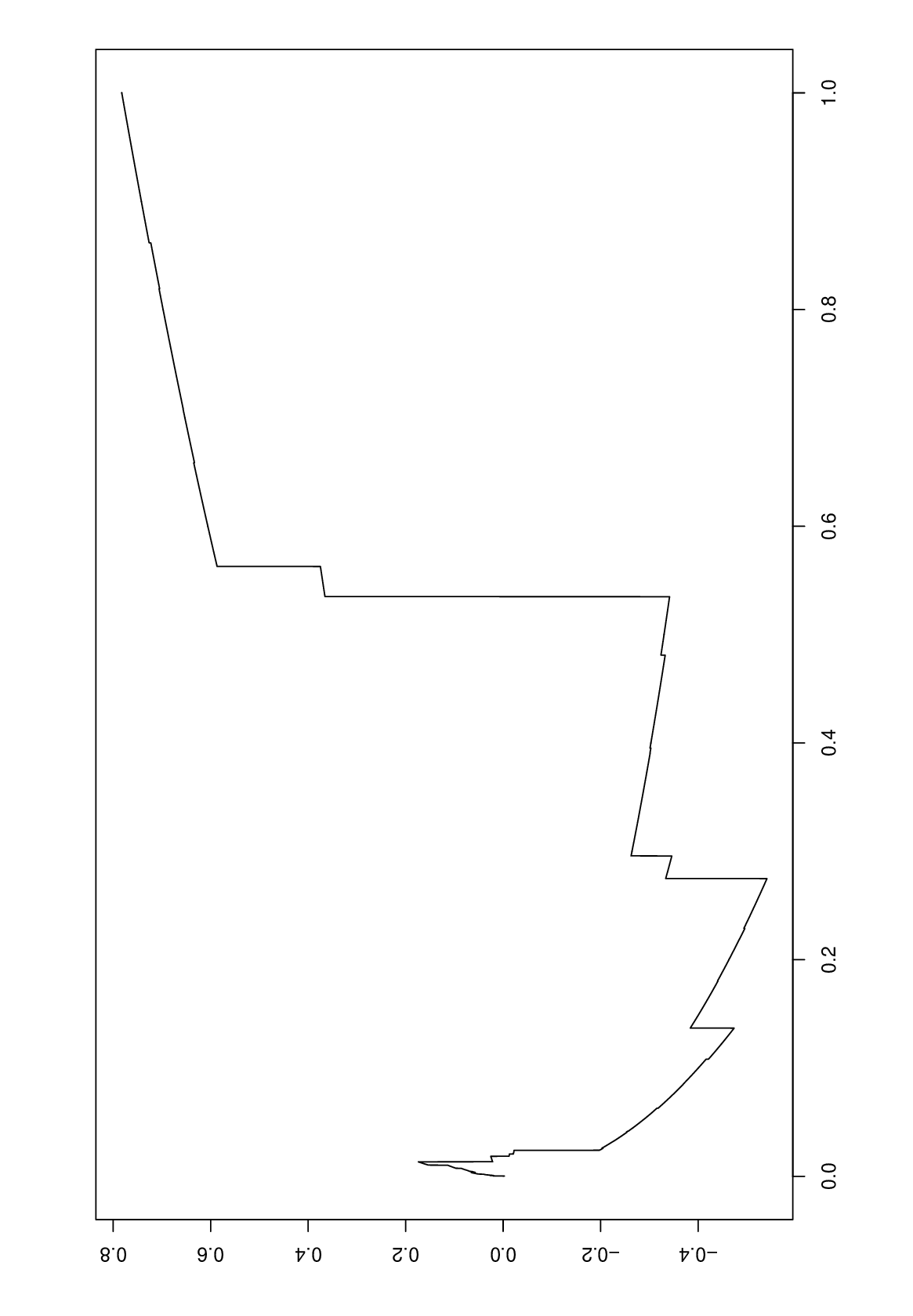}

\vspace*{6cm}

In this case it is possible to compute the generator for a much wider
class of functions.

\begin{proposition}
Let $G_\sigma$ be the log-convolution semi-group of the inverse log-gamma
distributions. Then (\ref{generator}) holds for any bounded function with
bounded first derivative.
\end{proposition}

{\bf Proof:} In the proof of Proposition \ref{generalcase}, we write that
$$\frac1{t-s}\left(\mathbb{E}[f(X_t)|X_s=x]-f(x)\right)
= \frac1s\frac1{\sigma^2-1}\int\big(\mathbb{E}\left[C_\sigma f(U_\sigma)\right]-f(x)\big)\phi(z)dz.$$
Denote by $\theta(u)$ the quantity $e^{-u/2}x+\sqrt s\sqrt{1-e^{-u}}z$. Then,
inserting $\mathbb{E}\left[C_1 f(U_\sigma)\right] =
\mathbb{E}\left[f\left(\theta(U_\sigma)\right)\right]$ we get,
\begin{eqnarray}
\lefteqn{\frac1{t-s}\left(\mathbb{E}[f(X_t)|X_s=x]-f(x)\right)}\label{gamma}\\
& = & \frac1s\frac1{\sigma+1}\int\left\{\frac{\mathbb{E}\left[C_\sigma f(U_\sigma)\right]
- \mathbb{E}\left[C_1 f(U_\sigma)\right]}{\sigma-1} +
\frac{\mathbb{E}\left[C_1 f(U_\sigma)\right] - f(x)}{\sigma-1}\right\}\phi(z)dz.
\nonumber
\end{eqnarray}
Since
$$\frac{C_\sigma f(U_\sigma) - C_1 f(U_\sigma)}{\sigma-1} =
\frac{f(\sigma\theta(U_\sigma)) - f(\theta(U_\sigma))}{\sigma-1} =
\theta(U_\sigma)f'(\eta_\sigma),$$
for some $\eta_\sigma$ between $\theta(U_\sigma)$ and $\sigma\theta(U_\sigma)$.
$\theta$ and $f'$ being bounded, we obtain that
$$\lim_{\sigma\downarrow1}\int\frac{\mathbb{E}\left[C_\sigma f(U_\sigma)\right]
- \mathbb{E}\left[C_1 f(U_\sigma)\right]}{\sigma-1}\phi(z)dz = xf'(x).$$
To compute the limit of the second term in (\ref{gamma}), we use Lemma \ref{gammalimit}
below, which shows that
\begin{eqnarray*}
\lefteqn{\lim_{\sigma\downarrow1}\int\frac{\mathbb{E}\left[C_1 f(U_\sigma)\right] - f(x)}{\sigma-1}\phi(z)dz}\\
& = &
a\int\int_0^\infty\frac{f\left(e^{-u/2}x+\sqrt s\sqrt{1-e^{-u}}z\right)-f(x)}ue^{-bu}du\phi(z)dz\\
& = &
\int[f(x+y)-f(x)]\int_0^\infty\phi(x(e^{-u/2}-1),s(1-e^{-u}),y)a\frac{e^{-bu}}ududy.
\end{eqnarray*}
\hfill$\Box$

Note that since $\nu((0,\infty))=+\infty$,
$\int_0^\infty\phi(x(e^{-u/2}-1),s(1-e^{-u}),y)\nu(du)du$ cannot be re-scaled to produce
a density for the jumps of the process.

\begin{lemma}
\label{gammalimit}
Let $V_p$ have a gamma distribution with density:
$$h_p(v)=\frac{b^p}{\Gamma(p)}v^{p-1}e^{-bv},\quad v>0.$$
Let $g$ be such that $g(0)=0$ and $g(v)/v$ is bounded. Then
$$\lim_{p\downarrow0}\frac{1}{p}\mathbb{E}[g(V_p)] =
\int_0^\infty\frac{g(v)}{v}e^{-bv}dv.$$
\end{lemma}

{\bf Proof:}
First observe that
$$\frac{1}{p}\mathbb{E}[g(V_p)] = \frac{1}{b}\mathbb{E}\left[\frac{g(V_{p+1})}{V_{p+1}}\right].$$
taking the limit as $p\downarrow0$, we obtain by dominated convergence
$$\lim_{p\downarrow0}\frac{1}{p}\mathbb{E}[g(V_p)] = \frac1b\mathbb{E}\left[\frac{g(V_1)}{V_1}\right]
= \int_0^\infty\frac{g(v)}{v}e^{-bv}dv.$$
\hfill$\Box$

{\bf Acknowledgement:} This research was supported by the Australian Research Council.
The authors would like to thank Boris Granovsky for very fruitful discussions during his
visit to Monash University.

\end{document}